\documentclass[12pt]{article}

\usepackage{latexsym}
\usepackage{amsmath}
\usepackage{graphicx}
\usepackage{subfigure}

\makeatletter
\newfont{\footsc}{cmcsc10 at 8truept}
\newfont{\footbf}{cmbx10 at 8truept}
\newfont{\footrm}{cmr10 at 10truept}
\makeatother
\title{Counting $d$-polytopes with $d+3$ vertices}

\author{\'Eric Fusy\\
\small Algorithm Project\\[-0.8ex]
\small INRIA Rocquencourt, France\\[-0.8ex]
\small \texttt{eric.fusy@inria.fr}}

\newtheorem{theorem}{Theorem}
\newtheorem{proposition}[theorem]{Proposition}
\newtheorem{lemma}[theorem]{Lemma}
\newtheorem{corollary}[theorem]{Corollary}

\begin{document}
\maketitle
\begin{abstract}
We completely solve the problem of enumerating combinatorially
inequivalent $d$-dimensional polytopes with $d+3$ vertices. A first
solution of this problem, by Lloyd, was published in 1970. But the
obtained counting formula was not correct, as pointed out in the new
edition of Gr\"unbaum's book. We both correct the mistake of Lloyd and
propose a more detailed and self-contained solution, relying on
similar preliminaries but using then a different enumeration method
involving automata. In addition, we introduce and solve the problem of
counting oriented and achiral (i.e. stable under reflection)
$d$-polytopes with $d+3$ vertices. The complexity of computing tables
of coefficients of a given size is then analyzed. Finally, we derive
precise asymptotic formulas for the numbers of $d$-polytopes, oriented
$d$-polytopes and achiral $d$-polytopes with $d+3$ vertices. This
refines a first asymptotic estimate given by Perles.
\end{abstract}
\emph{Acknowledgement.}
The author would like to thank G\"unter M. Ziegler for having posed the 
problem to him and having taken much time to explain thoroughly the 
combinatorics of Gale diagrams and to correct in details a first draft of 
the paper. 
\section{Introduction}
A \emph{polytope} $P$ is the convex hull of a finite set of point of a 
vector space $\mathbf{R}^d$. If $P$ is not contained in any hyperplane 
of $\mathbf{R}^d$, then $P$ is said $d$-dimensional, or is called a 
$d$-polytope. A \emph{vertex} (resp. a \emph{facet}) of $P$ is defined 
as the intersection of $P$ with an hyperplane $H$ of $\mathbf{R}^d$ such 
that $P\cap H$ has dimension $0$ (resp. has dimension $d-1$) and one of 
the two open sides of $H$ does not meet $P$. A vertex $v$ 
is \emph{incident} to a facet $f$ if $v\in f$.
 
This article addresses the problem of counting combinatorially 
different $d$-polytopes with $d+3$ 
vertices, meaning that two polytopes are identified if their incidences 
vertices-facets are isomorphic (i.e. the incidences are the same up to 
relabeling of the vertices). Whereas general $d$-polytopes are involved 
objects, $d$-polytopes with few vertices are combinatorially tractable. 
Precisely, each combinatorial $d$-polytope 
with $d+3$ vertices gives rise in a bijective way to a configuration of 
$d+3$ points in the plane, placed at the centre and at vertices of a 
regular $2k$-gon, and satisfying two local conditions and a global 
condition. 
As a consequence, counting combinatorial $d$-polytopes with $d+3$ vertices 
boils down to 
the much easier task of counting such configurations of $d+3$ points,
called reduced Gale diagrams. Following this approach,
Perles~\cite[p. 113]{grunbaum} gave
an explicit formula for the number of (combinatorial)   simplicial 
$d$-polytopes with
$d+3$ vertices and  Lloyd~\cite{LLoyd}
gave a more complicated formula for $c(d+3,d)$ the number of combinatorial 
$d$-polytopes with $d+3$ vertices. However, as pointed out in the
new edition of Gr\"unbaum's book~\cite[p. 121a]{grunbaum}, Lloyd's formula 
does not match 
with the first values of $c(d+3,d)$, obtained by
Perles~\cite[p. 424]{grunbaum}. 

In this article, we both correct the mistake of Lloyd and propose a more 
complete and self-contained solution for this enumeration problem. 
The following theorem is our main result:

\begin{theorem}
\label{theo:enum}
Let $c(d+3,d)$ be the number of combinatorially different $d$-polytopes 
with $d+3$ vertices. 
Then the generating function $\displaystyle P(x)=\sum_d c(d+3,d)x^{d+3}$ 
has the following 
expression, where $\phi(.)$ is the Euler totient function:
\begin{multline}
P(x)=-\frac{1}{1-x}\sum_{e\ odd}\frac{\phi(e)}{4e}\ln\left(1-\frac{2x^{3e}}{(1-2x)^{2e}}\right)+\frac{1}{1-x}\sum_{e\geq 1}\frac{\phi(e)}{2e}\ln\left(\frac{1-x^{e}}{1-2x^{e}}\right)\label{eq:compt}       \\
+{\frac {x \left( {x}^{2}-x-1 \right)  \left( {x}^{4}-{x}^{2}+1
 \right) }{ 2\left( 1-x \right)^2  \left( 2\,{x}^{6}-4\,{x}^{4}+4\,{x}^{
2}-1 \right) }}\nonumber-{\frac {x \left( {x}^{8}-2\,{x}^{7}+{x}^{6}+3\,{x}^{3}-{x}^{2}-x+1
 \right) }{ \left( 1+x \right) ^{2} \left( 1-x \right) ^{6}}}\nonumber
\end{multline}
\end{theorem}
The first terms of the series are $P(x)=x^5+7x^6+31x^7+116x^8+379x^9+1133x^{10}+3210x^{11}+\ldots$, i.e. there is one polytope with 5 vertices in the plane  
(the pentagon), there are 7 polytopes with 6 vertices in the 3-D space, etc.

The mistake of Lloyd, pointed precisely in
Section~\ref{sec:wheelnotP3}, is in the last rational term of
$P(x)$. Lloyd derived from his expression of $P(x)$ an explicit
formula for $c(d+3,d)$, which does not match with the correct values of
$c(d+3,d)$ because of the mistake in the computation of $P(x)$. We do not 
perform such a derivation for two
reasons: first, several equivalent formulas for $c(d+3,d)$ can be
derived from the expression of $P(x)$, so that the canonical form seems to be
on the generating function rather than on the
coefficients. Second, explicit formulas for $c(d+3,d)$ such as the one
of Lloyd involve double summations, hence require a quadratic number
of arithmetic operations to compute $c(d+3,d)$. As opposed to that and
discussed in Section~\ref{sec:complexity}, the coefficients $c(d+3,d)$
can be directly extracted iteratively from the expression of $P(x)$ in a 
very efficient
way: a table of the $N$ first coefficients can be computed with
$\mathcal{O}(N\log (N))$ operations. Using a mathematical
software like Maple, a table of several hundreds of coefficients can easily 
be
obtained.

In Section~\ref{sec:oriented}, we introduce the problem of counting 
\emph{oriented}
$d$-polytopes with $d+3$ vertices, meaning that two polytopes are
equivalent if they have the same combinatorial structure \emph{and} there exists an 
orientation-preserving of $\mathbf{R}^d$ homeomorphism mapping the first 
one to the second one. We establish a bijection
between oriented $d$-polytopes with $d+3$ vertices and
so-called oriented reduced Gale diagrams of size $d+3$, adapting the
original bijection so as to take the orientation into account. To our
knowledge, this oriented version of the bijection was not stated
before. The bijection implies that the task of counting oriented
$d+3$-vertex $d$-polytopes reduces to the task of counting oriented
reduced Gale diagrams with respect to the size, which is done in a
similar way as the enumeration of Gale diagrams. As a corollary, we
also enumerate combinatorial $d+3$-vertex $d$-polytopes giving rise
to only one oriented polytope. These polytopes, called \emph{achiral},
are also characterized as having a geometric representant fixed by a
reflection of $\mathbf{R}^d$.

Finally, in Section~\ref{sec:asym}, we give precise asymptotic estimates 
for the
coefficients $c(d+3,d)$, $c^{+}(d+3,d)$, $c^{-}(d+3,d)$ counting
(combinatorial) $d$-polytopes, oriented $d$-polytopes and achiral 
$d$-polytopes with
$d+3$ vertices. No asymptotic result was given in Lloyd's paper, but
Perles~\cite[p.114]{grunbaum} proved that there exist two constants $c_1$ 
and $c_2$ such
that $c_1\frac{\gamma^d}{d}\leq c(d+3,d)\leq c_2\frac{\gamma^d}{d}$,
where $\gamma$ is explicit, $\gamma\approx 2.83$. Using analytic combinatorics, we deduce from the
expression of $P(x)$ that $c(d+3,d)\sim c\frac{\gamma^d}{d}$, with $c$
an explicit constant and $\gamma$ equal to the $\gamma$ of Perles, but
with a simplified definition. Hence this agrees with Perle's estimate
and refines it.

\noindent\paragraph{Overview of the proof of Theorem~\ref{theo:enum}.} 
In Section~\ref{sec:def}, we  give a sketch of proof of the bijection 
between combinatorial $(d+3)$-vertex
$d$-polytopes and reduced Gale diagrams of size $d+3$. With this bijection,
 the enumeration of $(d+3)$-vertex $d$-polytopes reduces to the
enumeration of reduced Gale diagrams with respect to the size. 

The scheme of
our method of enumeration of reduced Gale diagrams follows, in a more
detailed way, the same lines as Lloyd. The first observation (see
Section~\ref{sec:reduced_gale_remarks}) is that it is sufficient to 
concentrate on the enumeration of reduced Gale diagrams with no label at 
the centre and satisfying 
the two local conditions (forgetting temporarily the third global
condition).  We introduce a special terminology for
these diagrams, calling them \emph{wheels}. As wheels are enumerated up 
to rotation and up to reflection, 
they are subject to symmetries: Burnside's lemma reduces the task of 
counting 
wheels to the task of counting so-called rooted wheels (where the presence 
of a root deletes possible symmetries) and rooted symmetric wheels of two 
types: rotation and reflection, see Section~\ref{sec:methwheels}. 

After these preliminaries, our
treatment for the enumeration of rooted wheels differs from that of
Lloyd, which relies on an auxiliary theorem of Read, requiring to
operate in two steps. The method we propose in 
Section~\ref{sec:countRooted} is direct and
self-contained: we associate with 
a rooted wheel a word on a specific (infinite) alphabet and we show that 
the set of words derived from rooted wheels is recognized by a simple 
automaton (see Figure~\ref{fig:generic_aut_a}). Under the framework of 
automata, generating functions appear as a very powerful tool providing 
simple 
(in general rational) and compact solutions in an automatic way. We derive 
from the automaton 
an explicit rational expression for the generating function of rooted wheels. 
The enumeration of rooted symmetric wheels is done in a similar way, 
associating words with such rooted wheels and observing that the obtained 
sets of words are recognized by automata. 
The injection into Burnside's Lemma of the rational expressions 
for rooted and rooted symmetric wheels yields an explicit expression for the generating functions of wheels, given in Section~\ref{sec:wheelsgen}. 
Theorem~\ref{theo:enum} follows after taking the global condition
(called half-plane condition) into account, which requires only some
exhaustive treatment of cases, see Section~\ref{sec:wheelnotP3}. 

\section{Gale diagrams of $(d+3)$-vertex $d$-polytopes}
\subsection{Gale diagrams}
\label{sec:def}
Following Perles and Lloyd, we define a \emph{reduced Gale diagram} as a 
regular 
$2k$-gon, with $k\geq 2$, that carries non-negative labels at its centre 
and at its vertices, 
with the following properties:

\begin{description}
\item[P1:] Two opposite vertices of the $2k$-gon cannot both have label 0
\item[P2:] Two neighbour vertices of the $2k$-gon cannot both have label 0
\item[P3:] (half-plane condition) Given any diameter of the $2k$-gon, the 
sum of 
the labels of vertices belonging to any (open) side of the diameter is at 
least 2.
\end{description}
In addition, two reduced Gale diagrams are identified if the first one can 
be 
obtained from the second one by a rotation or by a rotation and a reflection. 
The \emph{size} of a reduced Gale diagram is defined as the sum of its 
labels. The following theorem is essential in order to reduce the problem 
of enumeration of polytopes to the tractable problem of counting reduced 
Gale diagrams. Details of the proof can be found in Gr\"unbaum's 
book~\cite[Sect.~6.3]{grunbaum}. 

\begin{theorem}
\label{theo:bij}
\emph{(Perles)} The number of combinatorially different $d$-polytopes with $d+3$ 
vertices is equal to the number of reduced Gale diagrams of size $d+3$.
\end{theorem}
\noindent\emph{Proof} (Sketch): Given a $d$-polytope $P$ with $d+3$ 
vertices $v_1,\ldots,v_{d+3}$, a matrix $M_P$ is associated with $P$ in 
the following way: $M_P$ has $d+3$ rows, the $i$th row consisting of a 1 
followed by the $d$-vector-position of the vertex $v_i$. Hence, $M_P$ has 
$d+1$ columns, and it can be shown that $M_P$ has rank $d+1$. As a 
consequence, the vector space $\mathcal{V}(P)$ spanned by the column 
vectors $(C_1,\ldots,C_{d+1})$ 
of $M_P$ has dimension $d+1$, so its orthogonal 
$\mathcal{V}(P)^{\bot}$ has dimension 2. Let $(A_1,A_2)$ be a base of
$\mathcal{V}(P)^{\bot}$ and let $A$ be the $(d+3)\times 2$ 
matrix whose two columns are $(A_1,A_2)$. 
Then $A$ is called a Gale diagram of $P$. The matrix $A$ can be seen 
as a configuration of $d+3$ points in the plane, each point 
corresponding to a row of $A$. The combinatorial structure of $P$, 
i.e. the incidences vertices-facets, can be recovered from $A$. 
However, several Gale diagrams can correspond to the same (combinatorial) 
polytope. One can perform successive reductions, keeping the same 
associated combinatorial polytope, so that the $d+3$ points of 
the diagram are finally located either at the centre or at vertices 
of a regular $2k$-gon. Giving to the centre and to each vertex of the 
$2k$-gon a label indicating the number of points located at it, one 
obtains a $2k$-gon with labels characterized by the fact that they 
satisfy properties P1, P2 and P3. In addition, it can be shown that 
this reduction is maximal, i.e. that the combinatorial types of the 
polytopes 
associated with two inequivalent (i.e. not equal up to rotation and 
reflection) reduced Gale diagrams are different.
$\Box$

\subsection{Gale diagrams and wheels}
\label{sec:reduced_gale_remarks}
A first remark is that properties P1, P2, P3 do not depend on the value of 
the label at the centre of the $2k$-gon. Hence the number $g_n$ of reduced 
Gale diagrams of size $n$ is easily deduced from the coefficients $e_i$ 
counting reduced Gale diagrams of size $i$ with label 0 at the centre 
(such reduced Gale diagrams correspond to so-called non-pyramidal polytopes): 
$$
g_n=\sum_{i=1}^ne_i.
$$
As a consequence, we concentrate on the enumeration of labelled 
$2k$-gons   (meaning that only the $2k$ vertices of the $2k$-gon carry 
labels) 
satisfying properties P1, P2, P3.

A second remark is that Property P3 is implied by Property P2 if the number 
of diameters is at least 5. As a consequence, we will first put 
aside Property P3 and focus on the enumeration of labelled $2k$-gons 
satisfying properties P1 and P2 (counted up to rotation and up to 
reflections). Such labelled $2k$-gons are called \emph{wheels}. Wheels with 
2 vertices, even though corresponding to a degenerated polygon, 
are also counted. The enumeration of wheels will be performed in 
Section~\ref{sec:methwheels} and Section~\ref{sec:countwheels}. 
By definition of wheels, the number of reduced Gale diagrams with no label at the 
centre is obtained as the difference between the number of wheels and 
the number  of wheels not satisfying Property P3. The latter term, 
considered in Section~\ref{sec:wheelnotP3}, is easy to calculate 
using some exhaustive treatment of cases, because wheels not 
satisfying Property P3 have at most 4 diameters.

\section{Method of enumeration of wheels}
\label{sec:methwheels}  
\subsection{Rooted wheels}
\label{sec:rooted_wheels}

\begin{figure}
\begin{center}
\subfigure[A wheel.]{
\label{fig:wheel_a}
\includegraphics[scale=1.5]{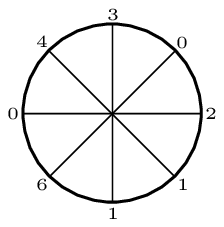}
}
\hspace{0.5cm}
\subfigure[A rooted wheel.]{
\label{fig:wheel_b}
\includegraphics[scale=1.5]{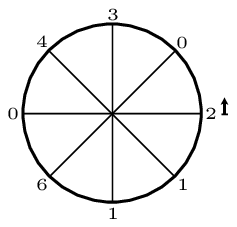}
}
\hspace{0.5cm}
\subfigure[The associated integer sequence.]{
\label{fig:wheel_c}
\mbox{$\ \ $\includegraphics[scale=1.5]{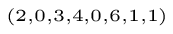}$\ \ $}
}
\caption{Example of wheel and rooted wheel.}
\label{fig:wheel}
\end{center}
\end{figure}

A wheel is \emph{rooted} by selecting one vertex of the $2k$-gon and by 
choosing a sense of traversal (clockwise or counter-clockwise) of the 
$2k$-gon. See Figure~\ref{fig:wheel_b} 
for an example~\footnote{On the figures, regular $2k$-gons 
are represented as $2k$ vertices regularly distributed on a circle, 
for aesthetic reasons and consistence with the terminology of wheels.}.

Traversing the $2k$-gon from the selected vertex in the direction 
indicated by the root, one obtains an integer sequence $(a_1,\ldots,a_{2k})$ 
satisfying the following conditions:

\begin{description}
\item[S1:] For each $1\leq i\leq 2k$, $a_i$ and $a_{(i+k)\mod 2k}$ are not 
both $0$.
\item[S2:] For each $1\leq i\leq 2k$, $a_i$ and $a_{(i+1)\mod 2k}$ are not 
both $0$.
\end{description}

An integer sequence satisfying properties S1 and S2 is called a 
\emph{wheel-sequence}. The \emph{size} of the wheel-sequence is defined 
as $(a_1+\ldots+a_{2k})$. Properties S1 and S2 are simply the 
respective transpositions of properties P1 and P2 on the integer sequence, 
so that we can identify rooted wheels with size $n$ and $k$ diameters and 
wheel-sequences of size $n$ and length $2k$.

\subsection{Burnside's lemma}
Burnside's lemma is a convenient tool to enumerate objects defined modulo 
the action of a group, which means that they are counted modulo symmetries. 
Let $G$ be a finite group acting on a finite set $E$. Given $g\in G$, 
we write $\mathrm{Fix}_g$ for the set of elements of $E$ fixed by $g$. 
Then the number of orbits of $E$ under the action of $G$ is given by:
\begin{equation}
|\mathrm{Orb}_E|=\frac{1}{|G|}\sum_{g\in G} |\mathrm{Fix}_g|
\end{equation}
where $|.|$ stands for cardinality.

\subsection{Burnside's lemma applied to wheels}
A wheel with size $n$ and $k$ diameters corresponds to an orbit of 
rooted wheels with size $n$ and $k$ diameters under the action of the 
dihedral group $D_{2k}$. Equivalently, using the identification between 
rooted wheels and wheel-sequences, a wheel with size $n$ and $k$ diameters 
corresponds to an orbit of wheel-sequences of size $n$ and length $2k$ under 
the action of $D_{2k}=\mathbf{Z}_{2k}\times \{+,-\}$, where the action 
is defined as follows, see Figure~\ref{fig:action}:

\begin{figure}
\begin{center}
\subfigure[Rotation action.]{
\label{fig:action_a}
\mbox{$\ \ \ \ \ \ \ \ \ \ \ \ \ \ \ $\includegraphics[scale=0.6]{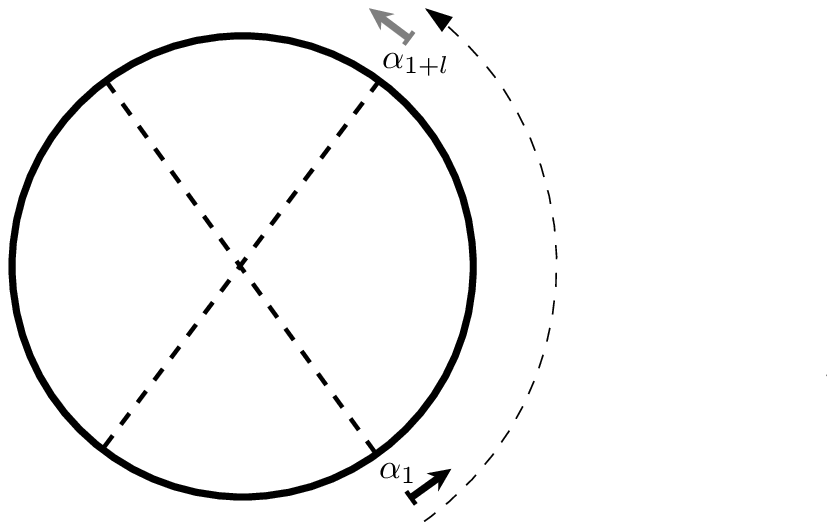}}
}
\subfigure[Reflexion action.]{
\label{fig:action_b}
\includegraphics[scale=0.6]{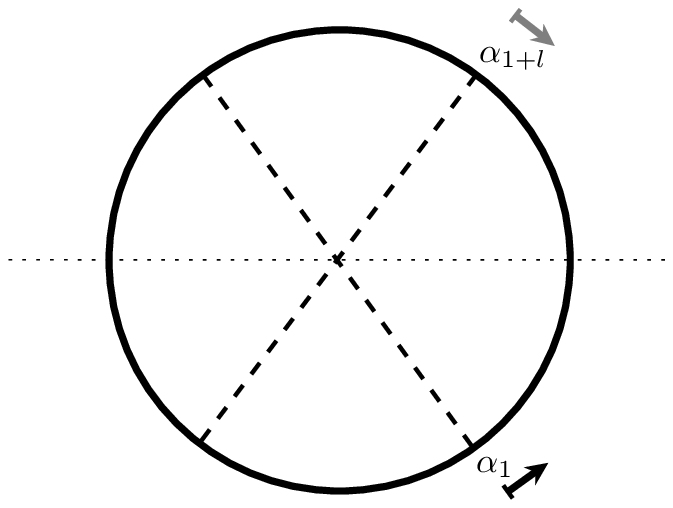}
}
\end{center}
\caption{The two cases of action of the dihedral group.}
\label{fig:action}
\end{figure}

\begin{eqnarray*}
(l,+)\cdot(a_1,\ldots,a_{2k})&=&(a_{1+l},\ldots,a_{2k},a_1,\ldots,a_l)\\
(l,-)\cdot(a_1,\ldots,a_{2k})&=&(a_{1+l},a_l,\ldots,a_1,a_{2k},\ldots,a_{2+l}),
\end{eqnarray*}
i.e. $(l,+)$ is a rotation and $(l,-)$ 
is a 
reflection.

Let us now introduce some definitions. A \emph{rotation-wheel} is a pair 
made of a rooted wheel and of a rotation of order at least 2 fixing the 
rooted wheel. Equivalently, it is a pair made of a sequence 
$(a_1,\ldots,a_{2k})$ and of an element $(l,+)$ with $l\neq 0$  
such that $(l,+)\cdot(a_1,\ldots,a_{2k})=(a_1,\ldots,a_{2k})$. 
A \emph{reflection-wheel} is a pair made of a rooted wheel and of 
a reflection fixing the rooted wheel. Equivalently, it is a pair made of a sequence $(a_1,\ldots,a_{2k})$ and of  an element $(l,-)$ 
such that  $(l,-)\cdot(a_1,\ldots,a_{2k})=(a_1,\ldots,a_{2k})$. The
following proposition ensures that, using Burnside's formula, counting
wheels reduces to counting rooted wheels, rotation wheels and
reflection wheels.

\begin{proposition}
\label{prop:burn}
Let $W_{n,k}$, $R_{n,k}$, $R_{n,k}^{+}$, $R_{n,k}^{-}$ be respectively 
the numbers of wheels, rooted wheels, rotation-wheels, and reflection-wheels 
with size $n$ and $k$ diameters. Let $W(x,u)$, $R(x,u)$, $R^{+}(x,u)$, and 
$R^{-}(x,u)$ be their generating functions. 
Then
\begin{equation}
\label{eq:burn}
4u\frac{\partial W}{\partial u}(x,u)=R(x,u)+R^{+}(x,u)+R^{-}(x,u).
\end{equation}
\end{proposition}
\noindent\emph{Proof:}
 As wheels with $k$ diameters are orbits 
of rooted wheels with $k$ diameters under the action of the dihedral 
group $D_{2k}$ (which has cardinality $4k$), Burnside's formula yields
$$
W_{n,k}=\frac{1}{4k}\left( R_{n,k}+R_{n,k}^{+}+R_{n,k}^{-}\right).
$$
Hence $\displaystyle \sum 4kW_{n,k}x^nu^k=\sum R_{n,k}x^nu^k+\sum R_{n,k}^{+}x^nu^k+\sum R_{n,k}^{-}x^nu^k$, which yields~(\ref{eq:burn}).
$\Box$

\section{Enumeration of wheels}
\label{sec:countwheels}
\subsection{Enumeration of rooted wheels}
\label{sec:countRooted}
In this section, we explain how to obtain a rational expression for the 
generating function $R(x,u)$ counting rooted wheels with respect to the 
size and number of diameters.

\subsubsection{The word associated to a rooted wheel.}
Let $s=(a_1,\ldots,a_{2k})$ be a wheel-sequence of size $n$ and length $2k$. 
Associate with $s$ the following word:

$$
\sigma:=\binom{a_1}{a_{k+1}},\binom{a_2}{a_{k+2}},\ldots,\binom{a_k}{a_{2k}}
$$

Observe that the length of $\sigma$ is the number of diameters of 
the associated rooted wheel. As each letter of $\sigma$ contains a pair of 
opposite vertices of the $2k$-gon, the fact that two opposite vertices are 
not both 0 (Property P1 or equivalently Property S1) translates into
the following property: 
$$\sigma \mathrm{\ is\ a\ word\ on\ the\ alphabet}\  
\mathcal{A}:=\mathbf{N}^2\backslash \left\{\binom{0}{0}\right\}.$$ 
Now let us detail the translation of Property  P2 (or S2) on the 
word $\sigma$. First, the alphabet $\mathcal{A}$ is partitioned into 
three subalphabets:

$$
\mathcal{B}=\left\{ \binom{i}{j}\ \mathrm{with}\ i>0,j>0\right\}, \mathcal{C}=\left\{ \binom{i}{0}\ \mathrm{with}\ i>0\right\}, \mathcal{D}=\left\{ \binom{0}{j}\ \mathrm{with}\ j>0\right\}
$$
Property S2 is translated as follows:
\begin{itemize}
\item
$a_i$ and $a_{i+1}$ are not both $0$ for $1\leq i\leq k-1$ 
$\Longleftrightarrow$ $\sigma_i$ and $\sigma_{i+1}$ 
are not both in $\mathcal{D}$ for $1\leq i\leq k-1$
\item
$a_{k+i}$ and $a_{k+i+1}$ are not both $0$ for $1\leq i\leq k-1$ 
$\Longleftrightarrow$ $\sigma_i$ and $\sigma_{i+1}$ are not both 
in $\mathcal{C}$ for $1\leq i\leq k-1$
\item
$a_k$ and $a_{k+1}$ are not both $0$ $\Longleftrightarrow$ the pair 
$(\sigma_1,\sigma_k)$ is not in $\mathcal{C}\times\mathcal{D}$
\item
$a_1$ and $a_{2k}$ are not both $0$  $\Longleftrightarrow$ the pair 
$(\sigma_1,\sigma_k)$ is not in $\mathcal{D}\times\mathcal{C}$
\end{itemize}

Hence $\sigma$ is characterized as a word on the alphabet $\mathcal{A}$ 
that contains no factor $\mathcal{C}\mathcal{C}$ nor factor 
$\mathcal{D}\mathcal{D}$ and such that the pair made of its first and 
last letter in not in $\mathcal{C}\times\mathcal{D}$ nor in 
$\mathcal{D}\times\mathcal{C}$.

The \emph{size} of a letter is defined as the sum of its two integers, 
and the size of the word $\sigma$ is defined as the sum of the sizes of 
its letters. Hence the size of a rooted wheel is equal to the size of 
its associated word.

Notice that the generating functions of the 
three subalphabets $\mathcal{B}$, $\mathcal{C}$, and $\mathcal{D}$ with 
respect to the size are:

\begin{equation}
\label{eq:subalpha}
B(x)=\left(\frac{x}{1-x}\right)^2\ \ \  C(x)=\frac{x}{1-x}\ \ \ D(x)=\frac{x}{1-x}
\end{equation}

\subsubsection{Generic automaton and its generating functions}
\begin{figure}
\begin{center}
\subfigure[The generic automaton.]{
\label{fig:generic_aut_a}
\includegraphics[width=5cm]{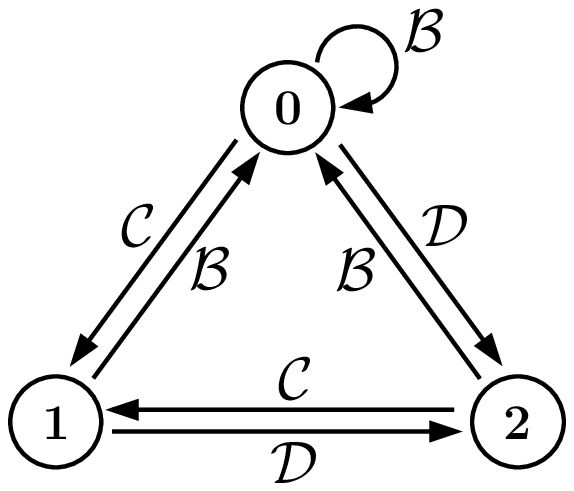}
}
\subfigure[Automaton recognizing words not containing $\mathcal{C}\mathcal{C}$ or $\mathcal{D}\mathcal{D}$.]{
\label{fig:generic_aut_b}
\mbox{$\ \ $\includegraphics[width=5cm]{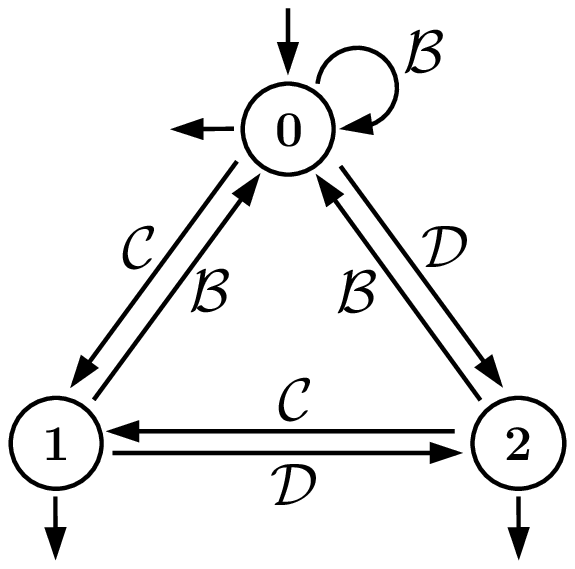}$\ \ $}
}
\end{center}
\caption{Automata associated with words not containing $\mathcal{C}\mathcal{C}$ or $\mathcal{D}\mathcal{D}$.}
\end{figure}

First we explain how to enumerate the words on the alphabet 
$\mathcal{A}$ avoiding the factors $\mathcal{C}\mathcal{C}$ and 
$\mathcal{D}\mathcal{D}$. The set of these words is 
recognized by the automaton represented on Figure~\ref{fig:generic_aut_b}, 
obtained from the automaton of Figure~\ref{fig:generic_aut_a} by 
choosing $\{0\}$ as starting state (entering arrow) and $\{0,1,2\}$ 
as end-states (leaving arrows). We call the automaton of 
Figure~\ref{fig:generic_aut_a} \emph{generic} because rooted wheels, 
rotation-wheels and reflection-wheels will give rise to languages 
on $\mathcal{A}$ recognized by slight modifications of this automaton.

For $i\in \{0,1,2\}$ and $j\in\{0,1,2\}$, we denote by 
$\mathcal{L}_{ij}$ the set of words accepted by the generic 
automaton that start at state $i$ and end at state $j$. 
Let $L_{ij}(x,u)$ be the generating function of $\mathcal{L}_{ij}$ 
with respect to the size and length of the word. Looking at the starting 
state and first letter of a word recognized by the generic automaton and 
ending at $0$, one gets the following system satisfied by the three 
generating functions $L_{00}(x,u)$, $L_{10}(x,u)$ and $L_{20}(x,u)$:

$$
\displaystyle
\left\{
\begin{array}{rcl}
L_{00}(x,u)&=&1+uB(x)L_{00}(x,u)+uC(x)L_{10}(x,u)+uD(x)L_{20}(x,u)\\
L_{10}(x,u)&=&uB(x)L_{00}(x,u)+uD(x)L_{20}(x,u)\\
L_{20}(x,u)&=&uB(x)L_{00}(x,u)+uC(x)L_{10}(x,u)
\end{array}
\right.
$$
Replacing $B(x)$, $C(x)$ and $D(x)$ by their expressions given 
in~(\ref{eq:subalpha}), this system becomes
$$
\left(
\begin{array}{c}
L_{00}(x,u)\\
L_{10}(x,u)\\
L_{20}(x,u)
\end{array}
\right)
=
\left(
\begin{array}{ccc}
\frac{ux^2}{(1-x)^2}&\frac{ux}{1-x}&\frac{ux}{1-x}\\
\frac{ux^2}{(1-x)^2}&0&\frac{ux}{1-x}\\
\frac{ux^2}{(1-x)^2}&\frac{ux}{1-x}&0
\end{array}
\right)\cdot
\left(
\begin{array}{c}
L_{00}(x,u)\\
L_{10}(x,u)\\
L_{20}(x,u)
\end{array}
\right)
+
\left(
\begin{array}{c}
1\\
0\\
0
\end{array}
\right).
$$
Solving this matrix equation, one gets explicit rational expressions 
for $L_{00}(x,u)$, $L_{10}(x,u)$ and $L_{20}(x,u)$, for instance: 
$$\displaystyle L_{10}(x,u)=\frac{ux^2(1-x)}{1-x(3+u-3x-ux+x^2+u^2x^2)}.$$ 
One can 
similarly define a matrix-equation satisfied by 
$\{ L_{01}(x,u),L_{11}(x,u),L_{21}(x,u)\}$ and a 
matrix-equation satisfied by  $\{ L_{02}(x,u),L_{12}(x,u),L_{22}(x,u)\}$, 
from which one gets explicit rational expressions for these generating 
functions.

\subsubsection{Expression of the generating function of rooted wheels}
As we have seen in Section~\ref{sec:rooted_wheels}, rooted wheels with size 
$n$ and $k$ diameters can be identified with non-empty words of size $n$ and length $k$ 
on the alphabet $\mathcal{A}$, avoiding the factors $\mathcal{C}\mathcal{C}$ 
and $\mathcal{D}\mathcal{D}$ and such that the pair made of their 
first and last letter is not in $\mathcal{C}\times\mathcal{D}$ nor in 
$\mathcal{D}\times\mathcal{C}$.

\begin{figure}
\begin{center}
\includegraphics[width=12cm]{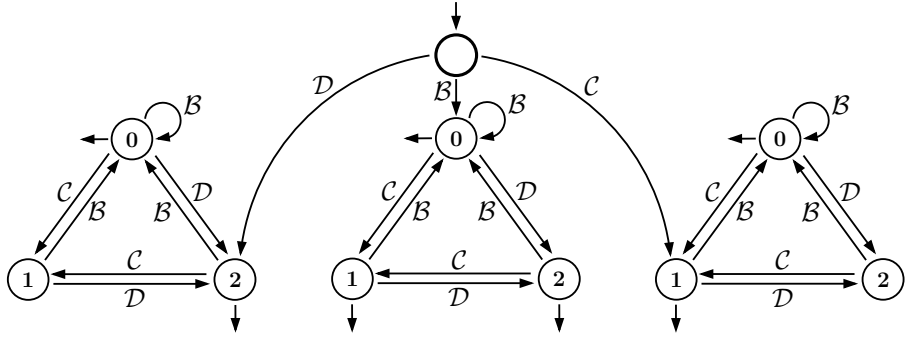}
\end{center}
\caption{Automaton recognizing non-empty words not containing 
$\mathcal{C}\mathcal{C}$ or $\mathcal{D}\mathcal{D}$ and not ending 
with $\mathcal{C}$ (resp. $\mathcal{D}$) if they start with $\mathcal{D}$ 
(resp. $\mathcal{C}$).}
\label{fig:aut_root}
\end{figure}

The language of these words is recognized by the automaton represented on 
Figure~\ref{fig:aut_root}. Hence the generating function $R(x,u)$ counting 
rooted wheels with respect to the size and number of diameters satisfies:

\begin{eqnarray*}
R(x,u)&=&uD(x)(L_{20}(x,u)+L_{22}(x,u))+uC(x)(L_{11}(x,u)+L_{10}(x,u))\\
&&+uB(x)(L_{00}(x,u)+L_{01}(x,u)+L_{02}(x,u))
\end{eqnarray*}
Replacing the generating functions on the right hand side by their rational 
expressions yields
\begin{equation}
\label{eq:exprrooted}
R(x,u)={\frac {ux \left( {u}^{2}{x}^{3}-2\,u{x}^{3}+2\,u{x}^{2}-{x}^{3}+4\,{
x}^{2}-5\,x+2 \right) }{ \left( {u}^{2}{x}^{3}-u{x}^{2}-3\,{x}^{2}+{x}
^{3}+3\,x+ux-1 \right)  \left( x-ux-1 \right) }}.
\end{equation}

\subsection{Enumeration of rotation-wheels}

As follows from the definition of rotation-wheels and from the 
identification 
between rooted wheels and wheel-sequences, a rotation-wheel corresponds 
to a pair made of a wheel-sequence $s=(a_1,\ldots,a_{2k})$ and of an 
element $l\in \mathbf{Z}_{2k}\backslash \{0\}$ such that the sequence $s$ is 
equal to its $l$-shift. Writing $e$ for the order of $l$ in 
$\mathbf{Z}_{2k}$ (hence $e$ divides $2k$), the sequence $s$ is 
characterized 
by the property that it can be written as $e$ concatenated copies of an 
integer sequence $(\alpha_1,\ldots,\alpha_{2k/e})$. In addition, it is 
well-known that for each divisor $e$ of $2k$ there are exactly $\phi(e)$ 
elements $l$ of order $e$ in $\mathbf{Z}_{2k}$. This yields the following 
lemma:

\begin{lemma}
\label{lemma:rot}
Let $\mathcal{R}^{(e)}$ be the set of rooted wheels whose wheel-sequence can 
be written as $e$ concatenated copies of an integer-sequence. Let 
$R^{(e)}(x,u)$  be the generating function of $\mathcal{R}^{(e)}$ with 
respect to the size and number of diameters. Then the generating function 
$R^{+}(x,u)$ of rotation-wheels is:
\begin{equation}
\label{eq:R+}
R^{+}(x,u)=\sum_{e\geq 2}\phi(e)R^{(e)}(x,u)
\end{equation}
\end{lemma} 

Let $e\geq 2$ and consider a rooted wheel of $\mathcal{R}^{(e)}$, so that 
its  associated sequence $(a_1,\ldots,a_{2k})$ consists of $e$ concatenated 
copies of an integer sequence $\alpha=(\alpha_1,\ldots,\alpha_{2k/e})$. 
We give a 
combinatorial characterization of the sequence  
$\alpha$ by distinguishing two cases:

\begin{figure}
\begin{center}
\subfigure[A rooted wheel of $\mathcal{R}^{(4)}$.]{
\label{fig:rotation_wheel_a}
\mbox{$\ \ \ \ $\includegraphics[scale=0.6]{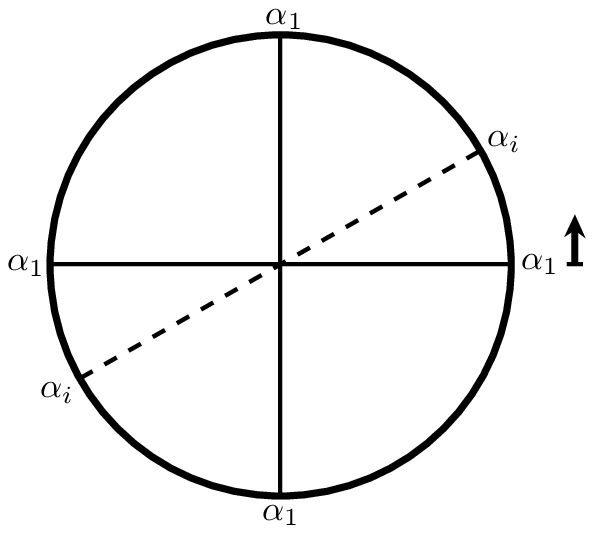}$\ \ \ \ $}
}
\hspace{1cm}
\subfigure[A rooted wheel of $\mathcal{R}^{(3)}$.]{
\label{fig:rotation_wheel_b}
\mbox{$\ \ \ \ $\includegraphics[scale=0.6]{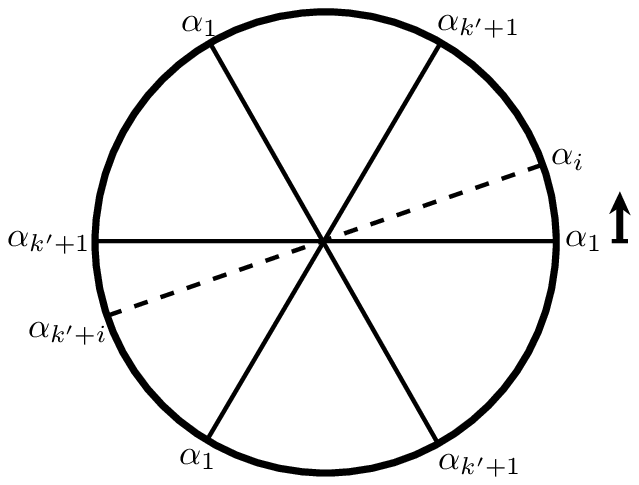}$\ \ \ \ $}
}
\end{center}
\caption{The two kinds of rooted wheels with a rotation-symmetry.}
\end{figure}

\noindent\paragraph{The number \boldmath{$e$} of copies is even.}
In this case, the opposite vertex of $\alpha_i$ on the $2k$-gon is 
$\alpha_i$, see Figure~\ref{fig:rotation_wheel_a}. As two opposite vertices 
of a wheel can not both have label 0 (Property P1), all integers $\alpha_i$ 
have to be positive. This condition ensures that two neighbour 
vertices of the $2k$-gon are not both 0 (Property P2).  Hence, for $r\geq 1$, 
a rooted wheel of $\mathcal{R}^{(2r)}$ with size $n$ and $k$ diameters 
corresponds to $2r$ concatenated copies of a non-empty sequence of positive 
integers of size $n/(2r)$ and length $k/r$. The generating function 
$I(x,u)$ counting non-empty sequences of positive integers with respect 
to the size and length is:
$$
I(x,u)=u\frac{x}{1-x}\frac{1}{1-u\frac{x}{1-x}}=\frac{ux}{1-x(1+u)}.
$$
Hence we obtain:
\begin{equation}
\label{eq:Reven}
R^{(2r)}(x,u)=I(x^{2r},u^r)=\frac{u^rx^{2r}}{1-x^{2r}(1+u^r)}.
\end{equation}

\noindent\paragraph{The number \boldmath{$e$} of copies is odd.} 
As $e$ is odd and divides $2k$, it also divides $k$. Hence $k/e$ is an 
integer, that we denote by $k'$. In this case, for $1\leq i\leq 2k'$, 
the opposite vertex of $\alpha_i$ on the $2k$-gon is 
$\alpha_{(i+k')\mod 2k'}$, see Figure~\ref{fig:rotation_wheel_b}. 
In addition, for $1\leq i\leq 2k'$, the next neighbour of $\alpha_i$ on 
the $2k$-gon is $\alpha_{(i+1)\mod 2k'}$.  As a consequence, the fact 
that $(a_1,\ldots,a_{2k})$ is a wheel-sequence is equivalent to the fact 
that $(\alpha_1,\ldots,\alpha_{2k'})$ is a wheel-sequence. Thus a 
wheel-sequence associated with a rooted wheel of $\mathcal{R}^{(2r+1)}$ 
corresponds to $(2r+1)$ concatenated copies of a wheel-sequence, 
so that for $r\geq 1$:
\begin{equation}
\label{eq:Rodd}
R^{(2r+1)}(x,u)=R(x^{2r+1},u^{2r+1}),
\end{equation}
where $R(x,u)$ is the generating function of rooted wheels.

Finally, equations~(\ref{eq:R+}),~(\ref{eq:Reven}) and~(\ref{eq:Rodd}) 
yield the following explicit expression 
of the generating function of rotation-wheels:
\begin{equation}
\label{eq:expR+}
R^{+}(x,u)=\sum_{r\geq 1}\phi(2r+1)R(x^{2r+1},u^{2r+1})+\sum_{r\geq 1}\phi(2r)\frac{u^rx^{2r}}{1-x^{2r}(1+u^r)}
\end{equation}

\subsection{Enumeration of reflection-wheels}
\label{sec:refl}
We recall that a reflection-wheel is a pair made of a rooted wheel and of 
a reflection fixing it. It can also be seen as a pair made of a 
wheel-sequence $(a_1,\ldots,a_{2k})$ and of an element $l\in \mathbf{Z}_{2k}$ 
such that $(a_1,\ldots,a_{2k})=(a_{1+l},a_l,\ldots,a_1,a_{2k},\ldots,a_{2+l})$

\begin{lemma}
\label{lemma:reflection_canonically}
Let $R^{(-1,-)}(x,u)$ be the generating function of rooted wheels fixed 
by the reflection $(-1,-)$ and let $R^{(0,-)}(x,u)$ be the generating 
function of rooted wheels fixed by  the reflection $(0,-)$. 
Then the generating function $R^{-}(x,u)$ of reflection-wheels is:
\begin{equation}
\label{eq:refl}
R^{-}(x,u)=u\frac{\partial}{\partial u}\left( R^{(-1,-)}(x,u)+R^{(0,-)}(x,u)\right).
\end{equation}
\end{lemma}
\noindent\emph{Proof:}
For $k\geq 1$ and $l\in \mathbf{Z}_{2k}$, 
we denote by $\mathcal{R}_{n,k}^{(l,-)}$ the set of rooted wheels 
with size $n$ and $k$ diameters whose associated sequence verifies 
$(a_1,\ldots,a_{2k})=(a_{1+l},a_l,\ldots,a_1,a_{2k},\ldots,a_{2+l})$. 
By definition, the set $\mathcal{R}_{n,k}^{-}$ of reflection-wheels 
with size $n$ and $k$ diameters is given by $\mathcal{R}_{n,k}^{-}=\cup_{l=0}^{2k-1}\mathcal{R}_{n,k}^{(l,-)}$. Observe that if a wheel sequence is 
fixed by the action of $(l,-)$, then its $r$-shift is fixed by the 
action of $(l-2r,-)$. As a consequence, $\mathcal{R}_{n,k}^{(l,-)}$ 
is in bijection with $\mathcal{R}_{n,k}^{(0,-)}$ if $l$ is 
even (these cases are those of a reflection fixing two vertices of 
the $2k$-gon); and $\mathcal{R}_{n,k}^{(l,-)}$ is in bijection 
with $\mathcal{R}_{n,k}^{(-1,-)}$ if $l$ is odd (these cases are 
those of a reflection fixing no vertex of the $2k$-gon). This directly 
yields $|\mathcal{R}_{n,k}^{-}|=k\left(|\mathcal{R}_{n,k}^{(0,-)}|+|\mathcal{R}_{n,k}^{(-1,-)}|\right)$, from which Equation~(\ref{eq:refl}) follows.
$\Box$

\subsubsection{Enumeration of rooted wheels fixed by the reflection $(0,-)$}
\label{sec:rooted_fixed0-}

\begin{figure}
\begin{center}
\subfigure[A rooted wheel fixed by the reflection $(0,-)$.]{
\label{fig:reflection_wheel_a}
\includegraphics[scale=0.6]{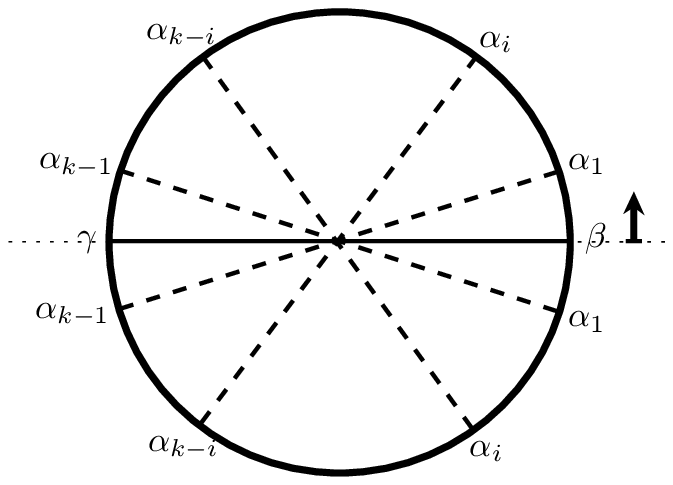}
}
\hspace{0.5cm}
\subfigure[A rooted wheel fixed by the reflection $(-1,-)$.]{
\label{fig:reflection_wheel_b}
\mbox{$\ $\includegraphics[scale=0.6]{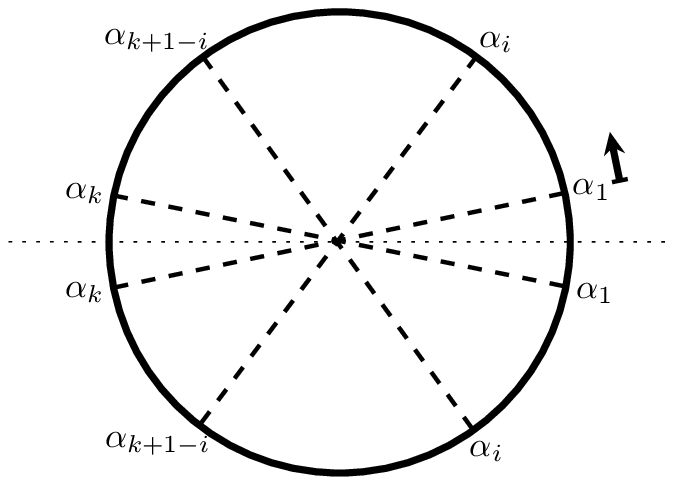}$\ $}
}

\end{center}
\caption{The two kinds of rooted wheels with a reflection-symmetry.} 
\end{figure}

Let $(a_1,\ldots,a_{2k})$ be a wheel-sequence fixed by $(0,-)$. 
Then $(a_1,\ldots,a_{2k})$ can be written as 
$(\beta,\alpha_1,\ldots,\alpha_{k-1},\gamma,\alpha_{k-1},\ldots,\alpha_1)$, 
see Figure~\ref{fig:reflection_wheel_a}. Observe that $\beta$ is opposite 
to $\gamma$ and that, for $1\leq i \leq  k-1$, $\alpha_i$ is opposite 
to $\alpha_{k-i}$ on the $2k$-gon. Then two cases arise:

\noindent\paragraph{The number of diameters is odd.} In this case, 
we write $r$ for the integer $(k-1)/2$. Then the fact that two opposite 
vertices of the rooted wheel do not both have label 0 (Property P1) is 
equivalent to the fact that:
$$
\sigma:=\binom{\beta}{\gamma},\binom{\alpha_1}{\alpha_{2r}},\ldots,\binom{\alpha_r}{\alpha_{r+1}}
$$ 
is a word on the alphabet $\mathcal{A}:=\mathbf{N}^2\backslash \{\binom{0}{0}\}$.

In addition, Property P2 (two neighbours can not both have label 0) 
translates as follows:

\begin{itemize} 
\item
$\beta$ and $\alpha_1$ are not both 0 and $\alpha_i$ and $\alpha_{i+1}$ 
are not both 0 for $1\leq i\leq r-1$ $\Longleftrightarrow$ $\sigma$ contains 
no factor $\mathcal{D}\mathcal{D}$.
\item
$\gamma$ and $\alpha_{2r}$ are not both 0 and $\alpha_i$ and $\alpha_{i+1}$ 
are not both 0 for $r+1\leq i\leq 2r-1$ $\Longleftrightarrow$ $\sigma$ 
contains 
no factor $\mathcal{C}\mathcal{C}$.
\item
$\alpha_r$ and $\alpha_{r+1}$ are not both 0 $\Longleftrightarrow$ the last 
letter of $\sigma$ is in the alphabet $\mathcal{A}$ (already implied by 
Property P1)
\end{itemize}

\begin{figure}
\begin{center}
\includegraphics[width=12cm]{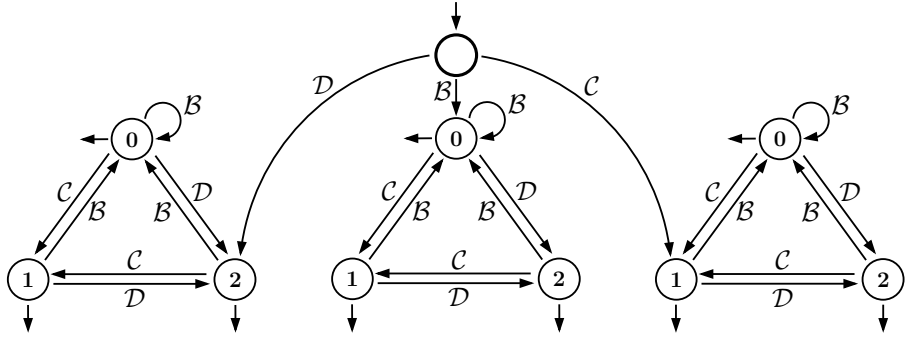}
\end{center}
\caption{Automaton recognizing words not containing $\mathcal{C}\mathcal{C}$ 
or $\mathcal{D}\mathcal{D}$, and reading the first letter separately.}
\label{fig:aut_refl}
\end{figure} 

Hence $\sigma$ is characterized as a non-empty word on the alphabet 
$\mathcal{A}$ that contains no factor $\mathcal{C}\mathcal{C}$ nor factor 
$\mathcal{D}\mathcal{D}$. The set of words on $\mathcal{A}$ satisfying this 
last property is already recognized by the automaton represented on 
Figure~\ref{fig:generic_aut_b}. However, the first letter of 
$\sigma$, containing the two fixed points, counts once in 
the sequence $(a_1,\ldots,a_{2k})$, whereas the other letters count twice. 
Hence we rather use the automaton of Figure~\ref{fig:aut_refl}, which
recognizes non-empty words avoiding $\mathcal{C}\mathcal{C}$ and 
$\mathcal{D}\mathcal{D}$, and where the first letter of the word is 
read separately. From this automaton, we get the generating function 
$R_{odd}^{(0,-)}(x,u)$ of rooted wheels fixed by $(0,-)$ and with an 
odd number of diameters:

\begin{eqnarray*}
R_{odd}^{(0,-)}(x,u)&=&uB(x)(L_{00}(x^2,u^2)+L_{01}(x^2,u^2)+L_{02}(x^2,u^2))\\
&&+uC(x)(L_{10}(x^2,u^2)+L_{11}(x^2,u^2)+L_{12}(x^2,u^2))\\
&&+uD(x)(L_{20}(x^2,u^2)+L_{21}(x^2,u^2)+L_{22}(x^2,u^2))\\
&=&{\frac { \left( x+1 \right)  \left( {x}^{3}+{u}^{2}{x}^{3}-2\,{x}^{2}-
x+2 \right) ux \left( 1-{x}^{2}+{x}^{2}u \right) }{ \left( x-1
 \right)  \left( {x}^{6}{u}^{4}-{u}^{2}{x}^{4}-3\,{x}^{4}+{x}^{6}+3\,{
x}^{2}+{u}^{2}{x}^{2}-1 \right) }}
\end{eqnarray*}

\noindent\paragraph{The number of diameters is even.}
The case of even $k$ is quite similar to the case of odd $k$. We denote 
by $r$ the integer $k/2$. Then, the fact that two opposite vertices of 
the $2k$-gon do not both have label 0 is equivalent to the fact that
$$
\sigma:=\binom{\beta}{\gamma},\binom{\alpha_1}{\alpha_{2r-1}},\ldots,\binom{\alpha_{r-1}}{\alpha_{r+1}}
$$ 
is a non-empty word on the alphabet $\mathcal{A}:=\mathbf{N}^2\backslash \{\binom{0}{0}\}$, and that $\alpha_r$ (which is self-opposite) is non 0. 
Similarly as for odd $k$, one can see that the $2k$-gon has no
neighbour vertices with label 0 iff the word $\sigma$ has no factor 
$\mathcal{C}\mathcal{C}$ nor factor $\mathcal{D}\mathcal{D}$. 
Hence the conditions for the word $\sigma$ (including the fact that the 
first letter of 
$\sigma$ is counted once and the other letters twice) are the same as 
for the words considered in the last paragraph, so that the generating 
function of these words is $R_{odd}^{(0,-)}(x,u)$. As a consequence, 
the generating function $R_{even}^{(0,-)}(x,u)$ of rooted wheels fixed by $(0,-)$ and with an even number of diameters verifies:
\begin{eqnarray*}
R_{even}^{(0,-)}(x,u)&=&R_{odd}^{(0,-)}(x,u)\frac{ux^2}{1-x^2}\\
&=&-{\frac { \left( x+1 \right) ^{2} \left( {x}^{3}+{u}^{2}{x}^{3}-2\,{x}
^{2}-x+2 \right) ux}{{x}^{6}{u}^{4}-{u}^{2}{x}^{4}-3\,{x}^{4}+{x}^{6}+
3\,{x}^{2}+{u}^{2}{x}^{2}-1}}
\end{eqnarray*}

Finally, the relation $R^{(0,-)}(x,u)=R_{odd}^{(0,-)}(x,u)+R_{even}^{(0,-)}(x,u)$ yields:
\begin{equation}
\label{eq:root0-}
R^{(0,-)}(x,u)={\frac {{u}^{2}{x}^{3} \left( x+1 \right)  \left( {x}^{3}+{u}^{2}{x}^{
3}-2\,{x}^{2}-x+2 \right) }{ \left( x-1 \right)  \left( {x}^{6}{u}^{4}
-{u}^{2}{x}^{4}-3\,{x}^{4}+{x}^{6}+3\,{x}^{2}+{u}^{2}{x}^{2}-1
 \right) }}
\end{equation}

\subsubsection{Enumeration of rooted wheels fixed by the reflection $(-1,-)$}
A wheel-sequence $(a_1,\ldots,a_{2k})$ fixed by $(-1,-)$ can be written 
as $(\alpha_1,\ldots,\alpha_k,\alpha_k,\ldots,\alpha_1)$. 
For $1\leq i\leq k-1$, the opposite vertex of $\alpha_i$ on the $2k$-gon 
is $\alpha_{k-i}$, see Figure~\ref{fig:reflection_wheel_b}. 
As in Section~\ref{sec:rooted_fixed0-}, two cases arise:

\noindent\paragraph{The number of diameters is odd.}
In this case, we write $r:=(k-1)/2$. As $\alpha_r$ is self-opposite 
and $\alpha_i$ is opposite to $\alpha_{2r+2-i}$ for $1\leq i\leq r-1$ 
the fact that two opposite vertices of the $2k$-gon have not both label 0 
is equivalent to the fact that:
$$
\sigma:=\binom{\alpha_1}{\alpha_{2r+1}},\ldots,\binom{\alpha_{r}}{\alpha_{r+2}}
$$
is a word (possibly empty) on the alphabet $\mathcal{A}$ and that 
$\alpha_{r+1}$ 
(which is self-opposite) is not 0. It is easily seen that the $2k$-gon 
has not two neighbour vertices with both label 0 iff the word $\sigma$ 
is empty or 
starts with a letter in $\mathcal{B}$ (because $\alpha_1$ and $\alpha_{k}$ 
are neighbour to themselves) and contains no factor $\mathcal{C}\mathcal{C}$ 
nor factor $\mathcal{D}\mathcal{D}$. Such words just consist of a letter 
in $\mathcal{B}$ followed by a word avoiding factors $\mathcal{C}\mathcal{C}$ 
and $\mathcal{D}\mathcal{D}$. As the set of words avoiding factors 
$\mathcal{C}\mathcal{C}$ 
and $\mathcal{D}\mathcal{D}$ is exactly recognized by the 
automaton of Figure~\ref{fig:generic_aut_b}, we can derive the following 
expression for the generating function $R_{odd}^{(-1,-)}(x,u)$ of rooted 
wheels fixed by $(-1,-)$ and with an odd  number of diameters:
\begin{eqnarray*}
R_{odd}^{(-1,-)}(x,u)&=&\left(1+u^2B(x^2)\cdot\left(L_{00}(x^2,u^2)+L_{01}(x^2,u^2)+L_{02}(x^2,u^2)\right)\right)\frac{ux^2}{1-x^2}\\
&=&-{\frac {{x}^{2}u \left( x-1 \right)  \left( x+1 \right)  \left( {x}^{
2}+{u}^{2}{x}^{2}-1 \right) }{{x}^{6}{u}^{4}-{u}^{2}{x}^{4}-3\,{x}^{4}
+{x}^{6}+3\,{x}^{2}+{u}^{2}{x}^{2}-1}}
\end{eqnarray*}

\noindent\paragraph{The number of diameters is even.}
In this case, we write $r:=k/2$. Then the $2k$-gon has no opposite vertices 
both carrying label 0 iff
$$
\sigma:=\binom{\alpha_1}{\alpha_{2r}},\ldots,\binom{\alpha_{r}}{\alpha_{r+1}}
$$
is a word on the alphabet $\mathcal{A}$. It is then easily seen that the 
$2k$-gon has no neighbour vertices both carrying label 0 iff $\sigma$ 
avoids the factors $\mathcal{C}\mathcal{C}$ and $\mathcal{D}\mathcal{D}$ 
and starts with a letter in $\mathcal{B}$. As mentioned above, such words 
consist of a letter in $\mathcal{B}$ followed by a word recognized by the 
automaton of Figure~\ref{fig:generic_aut_b}. Hence the generating function 
$R_{even}^{(-1,-)}(x,u)$ of rooted wheels fixed by $(-1,-)$ and with an 
even number of diameters is
\begin{eqnarray*}
R_{even}^{(-1,-)}(x,u)&=&u^2B(x^2)\left(L_{00}(x^2,u^2)+L_{01}(x^2,u^2)+L_{02}(x^2,u^2)\right)\\
&=&{\frac { \left( {x}^{2}-{u}^{2}{x}^{2}-1 \right) {u}^{2}{x}^{4}}{{u}
^{4}{x}^{6}-{u}^{2}{x}^{4}-3\,{x}^{4}+{x}^{6}+3\,{x}^{2}+{u}^{2}{x}^{2
}-1}}
\end{eqnarray*}

Finally, the relation $R^{(-1,-)}(x,u)=R_{odd}^{(-1,-)}(x,u)+R_{even}^{(-1,-)}(x,u)$ yields
\begin{equation}
R^{(-1,-)}(x,u)={\frac {u{x}^{2} \left( u{x}^{4}-{u}^{3}{x}^{4}-u{x}^{2}-1+2\,{x}^{2
}-{x}^{4}+{u}^{2}{x}^{2}-{u}^{2}{x}^{4} \right) }{{u}^{4}{x}^{6}-{u}^{
2}{x}^{4}-3\,{x}^{4}+{x}^{6}+3\,{x}^{2}+{u}^{2}{x}^{2}-1}}.
\end{equation}

\subsection{Expression of the generating function of wheels}
\label{sec:wheelsgen}
From Burnside's formula and from the expressions of the generating functions 
of rooted wheels, rotation-wheels and reflection-wheels, an explicit 
expression can be derived for the generating function of wheels:

\begin{proposition}
\label{prop:W}
Let $W(x)$ be the generating function of wheels with respect to the size. 
Then $W(x)$ has the following expression:
\begin{eqnarray}
W(x)&=&-\sum_{e\ odd}\frac{\phi(e)}{4e}\ln\left(1-\frac{2x^{3e}}{(1-2x)^{2e}}\right)+\sum_{e\geq 1}\frac{\phi(e)}{2e}\ln\left(\frac{1-x^{e}}{1-2x^{e}}\right)\nonumber\\&&+{\frac {x \left( {x}^{2}-x-1 \right)  \left( {x}^{4}-{x}^{2}+1
 \right) }{ 2\left( 1-x \right)  \left( 2\,{x}^{6}-4\,{x}^{4}+4\,{x}^{
2}-1 \right) }},\label{eq:exp_W}
\end{eqnarray}
where $\phi(.)$ is Euler totient function.
\end{proposition} 
\noindent\emph{Proof:} A first easy observation is that $W(x)=W(x,u)\big|_{u=1}$, where $W(x,u)$ is the generating function of wheels with respect to 
the size and number of diameters. Notice also that $W(x,0)=0$, because a wheel has at least one diameter, so $W(x)=\int_0^1\frac{\partial W}{\partial u}(x,u)du$. Hence the expression of 
$\frac{\partial W}{\partial u}(x,u)$ given in Proposition~\ref{prop:burn} yields
$$
W(x)=\frac{1}{4}\int_0^1 u^{-1}R(x,u)du+\frac{1}{4}\int_0^1 u^{-1}R^{+}(x,u)du+\frac{1}{4}\int u^{-1}R^{-}(x,u)du.
$$

According to Proposition~\ref{lemma:reflection_canonically}, $R^{-}(x,u)=u\frac{\partial}{\partial u}\left(R^{(-1,-)}(x,u)+R^{(0,-)}(x,u)\right)$. 
Hence, the explicit expressions of $R^{(0,-)}(x,u)$ and $R^{(-1,-)}(x,u)$ obtained in Section~\ref{sec:refl} yield
\begin{eqnarray*}
\frac{1}{4}\int_0^1 u^{-1}R^{-}(x,u)du&=&\frac{1}{4}\left(R^{(-1,-)}(x,1)+R^{(0,-)}(x,1)\right)\\
&=&{\frac {x \left( {x}^{2}-x-1 \right)  \left( {x}^{4}-{x}^{2}+1
 \right) }{ 2\left( 1-x \right)  \left( 2\,{x}^{6}-4\,{x}^{4}+4\,{x}^{
2}-1 \right) }}.
\end{eqnarray*}

Then, the expression of $R^{+}(x,u)$ given in~(\ref{eq:expR+}) yields:
\begin{eqnarray*}
\int_0^1 u^{-1}(R(x,u)+R^{+}(x,u))du&=&\sum_{r\geq 0}\phi(2r+1)\int_0^1 u^{-1}R\left(x^{2r+1},u^{2r+1}\right)du\\
&&+\sum_{r\geq 1}\phi(2r)\int_0^1 u^{-1}\frac{u^rx^{2r}}{1-x^{2r}-u^rx^{2r}}du.
\end{eqnarray*}
Writing $G(x):=\int_0^1 u^{-1}R(x,u)du$ and $H(x):=\int_0^1 u^{-1}\frac{ux^2}{1-x-ux^2}du$, and using the change of variable $y=u^{2r+1}$ for the 
$r$th term of the first sum and $y=u^r$ for the $r$th term of the 
second sum, we obtain
$$
\int_0^1 u^{-1}(R(x,u)+R^{+}(x,u))du=\sum_{r\geq 0}\frac{\phi(2r+1)}{2r+1}G\left(x^{2r+1}\right)+\sum_{r\geq 1}\frac{\phi(2r)}{r}H(x^r),
$$
so that
\begin{eqnarray*}
W(x)&=&\frac{1}{4}\left(\sum_{r\geq 0}\frac{\phi(2r+1)}{2r+1}G(x^{2r+1})+\sum_{r\geq 1}\frac{\phi(2r)}{r}H(x^r)\right)\\
&&+{\frac {x \left( {x}^{2}-x-1 \right)  \left( {x}^{4}-{x}^{2}+1
 \right) }{ 2\left( 1-x \right)  \left( 2\,{x}^{6}-4\,{x}^{4}+4\,{x}^{
2}-1 \right) }}
\end{eqnarray*}
The integral $H(x)$ is easy to compute: $H(x)=-\ln  \left( 1-2\,{x}^{2} \right) +\ln  \left( 1-{x}^{2} \right)$. Using Expression~(\ref{eq:exprrooted}) 
of $R(x,u)$ and a mathematical package, one finds $G(x)=-\ln  \left( 1-4\,x+4\,{x}^{2}-2\,{x}^{3} \right) +2\,\ln  \left( 1-x
 \right)$. Then, as observed by Lloyd, $\ln(1-4x+4x^2-2x^3)=2\ln(1-2x)+\ln(1-2x^3/(1-2x)^2)$, so that the terms of $W(x)$ can be re-combined into~(\ref{eq:exp_W}).
$\Box$

\section{Enumeration of wheels not satisfying P3}
\label{sec:wheelnotP3}
As observed in Section~\ref{sec:reduced_gale_remarks}, wheels not 
satisfying the half-plane condition P3 have at most four diameters. 
Hence the generating function of these wheels is equal to the difference 
between the generating function of wheels with at most four diameters and 
the generating function of wheels satisfying Property P3 and having at 
most four diameters.

\begin{figure}
\begin{center}
\includegraphics[width=12cm]{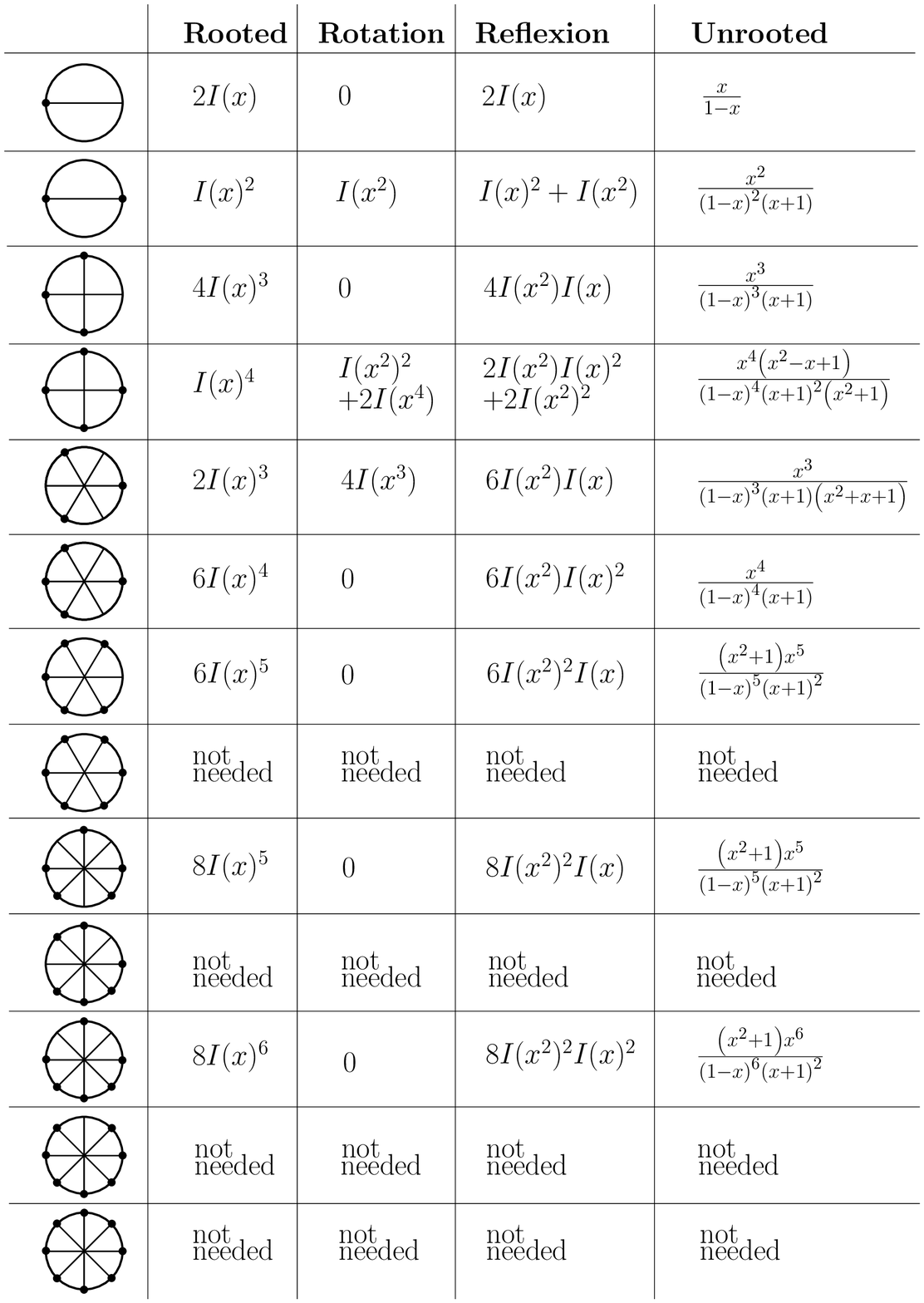}
\end{center}
\caption{The 13 possible configurations of a wheel with at most 4 diameters. 
A vertex has a disk iff its label is positive. Here $I(x):=x/(1-x)$ stands 
for the generating function of positive integers.}
\label{fig:13cas}
\end{figure}

\begin{figure}
\begin{center}
\includegraphics[width=12cm]{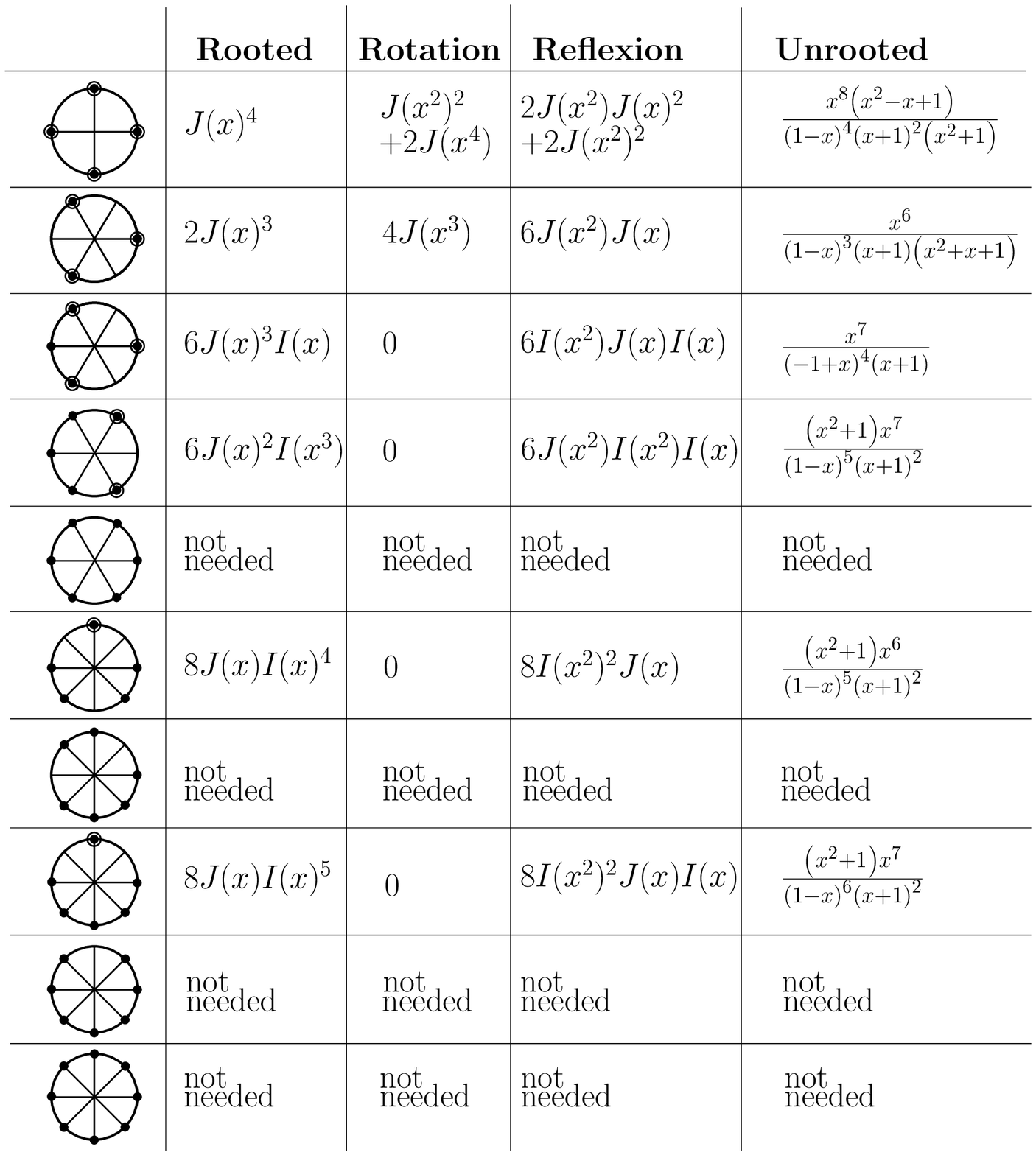}
\end{center}
\caption{The 10 possible configurations of a wheel satisfying the 
half-plane property and having at most 4 diameters. A vertex has no 
disk if it has label 0, has a black disk if its label is positive, 
and has a rounded black disk if its label must be at least 2 in order 
to satisfy P3. 
Here $I(x):=x/(1-x)$ stands for the generating function of positive 
integers and $J(x):=x^2/(1-x)$ stands for the generating function of 
integers of value at least 2.}
\label{fig:10cas}
\end{figure}

The \emph{configuration} of a wheel is obtained by putting a black disk on 
each vertex of the $2k$-gon occupied by a positive label and then by 
removing the labels.
Figure~\ref{fig:13cas} features the 13 possible configurations of a wheel 
with at most 4 diameters.  Similarly, Figure~\ref{fig:10cas} shows the 10 
possible configurations of a wheel with at most 4 diameters 
and satisfying P3, where a black disk is rounded if the label of 
the corresponding vertex must be at least 2 in order to satisfy 
Property P3. For each case on Figure~\ref{fig:13cas} and 
Figure~\ref{fig:10cas}, we can calculate the generating functions 
of rooted wheels, rotation-wheels and reflection-wheels corresponding 
to this configuration, and derive from Burnside's formula the generating 
function of wheels having this configuration. For example, the contribution 
of the 5th case of Figure~\ref{fig:13cas} is $\frac{1}{12}\left( 2I(x)^3+4I(x^3)+6I(x^2)I(x)\right)$, where $I(x)=x/(1-x)$. 

Then the generating function of wheels not satisfying Property P3 is 
obtained by taking the difference between the sum of the 13 contributions 
of Figure~\ref{fig:13cas} (last column) and the sum of the 10 contributions 
of Figure~\ref{fig:10cas} (last column). Observe that Cases 8, 10, 12, 13 of 
Figure~\ref{fig:13cas} exactly match Cases 5, 7, 9, 10 of 
Figure~\ref{fig:10cas}. Hence it is not necessary to compute the 
generating functions of these cases as they disappear in the difference. 
The calculation of the difference yields the following expression of the 
generating function $W_{\overline{P3}}(x)$ of wheels not satisfying 
Property P3:
\begin{equation}
\label{eq:wheel_not_P3}
W_{\overline{P3}}(x)={\frac {x \left( {x}^{8}-2\,{x}^{7}+{x}^{6}+3\,{x}^{3}-{x}^{2}-x+1
 \right) }{ \left( x+1 \right) ^{2} \left( 1-x \right) ^{5}}}.
\end{equation}

\noindent\paragraph{Remark:} Lloyd did a mistake in the calculation of $W_{\overline{P3}}(x)$. Precisely he forgot to subtract the term corresponding to the 8th 
configuration of Figure~\ref{fig:10cas} in his computation of the generating 
function of wheels not satisfying Property P3 and having 4 diameters and 
2 vertices with label 0. His presentation also has two typos: the label $\geq 2$ at the 
top of the top-right diagram of Fig.4 (page 129) has to be replaced by 
$\geq 1$ (it seems it is just a typo as the corresponding generating 
function is then correctly calculated). The second typo is in the term 
$(2,3)$ of page 131, where $\frac{1}{(1-x^4)}$ has to be replaced by 
$\frac{1}{(1-x)^4}$.

\section{Proof of Theorem~\ref{theo:enum}}
\label{sec:proof}
By definition (see Section~\ref{sec:def}), a reduced Gale diagram with no 
label at the centre is a wheel satisfying the half-plane property P3. 
Hence, the generating function of reduced Gale diagrams with no label at 
the centre is the difference between the generating function of wheels, 
given in~(\ref{eq:exp_W}), and the generating function of wheels not 
satisfying Property P3, given in~(\ref{eq:wheel_not_P3}). Then, as 
discussed in Section~\ref{sec:reduced_gale_remarks}, the generating 
function of reduced Gale diagrams is the multiplication by $1/(1-x)$ 
of the generating function of reduced Gale diagrams with no label at 
the centre. Finally, according to Theorem~\ref{theo:bij}, the number of 
reduced  Gale diagrams of size $d+3$ is equal to the number $c(d+3,d)$ of 
combinatorial $d$-polytopes with $d+3$ vertices. This yields 
Theorem~\ref{theo:enum}.

\section{Oriented and achiral \boldmath{$d$}-polytopes with \boldmath{$d+3$} vertices}
\label{sec:oriented}
This section deals with the enumeration of oriented $d$-polytopes with $d+3$ vertices, meaning that two polytopes $P$ and $P'$ are identified if they 
have the same combinatorial type \emph{and} there exists an 
orientation-preserving homeomorphism mapping $P$ to $P'$. We also introduce 
here \emph{oriented reduced Gale diagrams}, meaning that two reduced Gale 
diagrams are identified if they differ by a rotation.

\begin{theorem}
\label{theo:orie}
Oriented $d$-polytopes with $d+3$ vertices are in bijection with oriented 
reduced Gale diagrams of size $d+3$.
\end{theorem} 

\noindent\emph{Proof:(Sketch)}
The sketch of proof is very similar to the one of Theorem~\ref{theo:bij}. 
Hence
we keep the same notations, i.e. the matrix $M_P$ associated with $P$
and the vector space $\mathcal{V}(P)$ spanned by the column vectors
$(C_1,\ldots,C_{d+1})$ of $M_P$. An \emph{oriented} Gale diagram of
$P$ is a $(d+3)\times 2$ matrix whose two columns vectors $A_1$ and
$A_2$ form a base of $\mathcal{V}(P)^{\bot}$ and verify
$Det(C_1,\ldots,C_{d+1},A_1,A_2)>0$. As mentioned in the proof of
Theorem~\ref{theo:bij}, the combinatorial structure of $P$ is encoded in $A$, and also in
the reduced form of $A$. More precisely, if two polytopes have the
same combinatorial structure and the same orientation, then they have
the same reduced oriented Gale diagram. In addition, if two polytopes
$P$ and $P'$ have equivalent (i.e. equal up to \emph{rotation} only) 
oriented reduced Gale diagrams,
then one can deduce from it a continuous deformation of $P$ into $P'$,
keeping a polytope equivalent to $P$ all the way. In particular, $P$
and $P'$ have the same combinatorial structure and same orientation.
$\Box$

This bijection ensures that counting oriented $d$-polytopes with $d+3$
vertices reduces to counting oriented Gale diagrams with respect to the size.
This task is done in the same way as the enumeration of Gale diagrams,
i.e. we first enumerate oriented wheels and then substract oriented
wheels not satisfying the half-plane property P3. The only difference
between wheels and oriented wheels is in the application of Burnside's 
lemma. Namely, oriented wheels with $k$ diameters correspond to orbits of 
rooted wheels with $k$ diameters under the action of the cyclic group 
$\mathbf{Z}_{2k}$, of cardinality $2k$. From Burnside's lemma applied 
to the group $\mathbf{Z}_{2k}$, it follows that the generating 
function $W^{+}(x,u)$ of oriented wheels with respect to the size and 
number of diameters satisfies:
$$ 
2u\frac{\partial W^{+}}{\partial u}(x,u)=R(x,u)+R^{+}(x,u),
$$
where $R(x,u)$ and $R^{+}(x,u)$ are the generating functions of rooted 
wheels and of rotation-wheels.
Proceeding in a similar way as in the proof of Proposition~\ref{prop:W}, 
the following expression is obtained for the generating function $W^{+}(x)$ 
of oriented wheels with respect to the size:

$$
W^{+}(x)=-\sum_{e\ odd}\frac{\phi(e)}{2e}\ln\left(1-\frac{2x^{3e}}{(1-2x)^{2e}}\right)+\sum_{e\geq 1}\frac{\phi(e)}{e}\ln\left(\frac{1-x^{e}}{1-2x^{e}}\right)
$$

Then, oriented wheels not satisfying the half-plane property P3 are 
enumerated by doing the same exhaustive treatment of cases as in 
Section~\ref{sec:wheelnotP3}. For each of the 13 configurations of 
Figure~\ref{fig:13cas} and each of the 10 configurations of 
Figure~\ref{fig:10cas}, the associated generating function of oriented 
wheels is calculated using Burnside's Lemma (oriented formulation). 
For example the contribution of the second case of Figure~\ref{fig:10cas} 
is $\frac{1}{6}\left(2J(x)^3+4J(x^3)\right)={\frac { \left( {x}^{2}-x+1 \right) {x}^{6}}{ \left( 1-x \right) ^{3
}}}$. The generating function $W_{\overline{P3}}^{+}(x)$ of oriented wheels 
not satisfying P3 is the difference between the 13 oriented contributions 
of Figure~\ref{fig:13cas} and the 10 oriented contributions of 
Figure~\ref{fig:10cas}. This yields
$$W_{\overline{P3}}^{+}(x)={\frac { {x}^{11}+3\,{x}^{10}-3\,{x}^{9}-7\,{x}^{8}+4\,{x}^{7}+
4\,{x}^{6}+4\,{x}^{5}+3\,{x}^{4}-2\,{x}^{3}+x }{ \left( x+1
 \right) ^{3} \left( 1-x \right) ^{5}}}.
$$

Then, the generating function of oriented reduced Gale diagrams is equal to 
$$\frac{1}{1-x}\left(W^{+}(x)-W_{\overline{P3}}^{+}(x)\right),$$ 
see Section~\ref{sec:reduced_gale_remarks} and Section~\ref{sec:proof} 
for an explanation. As oriented Gale diagrams of size $d+3$ are in  
bijection with oriented $d$-polytopes with $d+3$ vertices, we obtain 
the following enumerative result:

\begin{theorem}
\label{theo:enum_oriented}
Let $c^{+}(d+3,d)$ be the number of oriented $d$-polytopes with $d+3$ 
vertices. 
Then the generating function $P^{+}(x):=\sum_dc^{+}(d+3,d)x^{d+3}$ 
has the following 
expression, where $\phi(.)$ is Euler totient function:
\begin{eqnarray*}
P^{+}(x)&=&\frac{1}{1-x}\left(-\sum_{e\ odd}\frac{\phi(e)}{2e}\ln\left(1-\frac{2x^{3e}}{(1-2x)^{2e}}\right)+\sum_{e\geq 1}\frac{\phi(e)}{e}\ln\left(\frac{1-x^{e}}{1-2x^{e}}\right)\right)\nonumber\\
&&-{\frac { {x}^{11}+3\,{x}^{10}-3\,{x}^{9}-7\,{x}^{8}+4\,{x}^{7}+
4\,{x}^{6}+4\,{x}^{5}+3\,{x}^{4}-2\,{x}^{3}+x }{ \left( x+1
 \right) ^{3} \left( 1-x \right) ^{6}}}.
\label{eq:exp_reduced}
\end{eqnarray*}
\end{theorem}
 The first terms of the series are $P^{+}(x)=x^5+7x^6+38x^7+170x^8+617x^9+1979x^{10}+5859x^{11}+\ldots$.

Observe 
that a $d$-polytope with $d+3$ vertices either gives rise to two different oriented
polytopes or to one oriented polytope. In the first (resp. second)
case, the polytope is called \emph{chiral} (resp. \emph{achiral}). It
can be shown that a combinatorial $(d+3)$-vertex $d$-polytope is
achiral iff one of its geometric representations is invariant under
the reflection $x_1\to -x_1$ (the proof relies on the fact that
achiral $d+3$-vertex $d$-polytopes are in bijection with reduced Gale
diagram having a reflection-symmetry). It follows from the definition that the generating function of achiral  
$d$-polytopes with $d+3$ vertices is equal to $2P(x)-P^{+}(x)$ where $P(x)$ 
and $P^{+}(x)$ are respectively the generating function of $d$-polytopes 
and of oriented $d$-polytopes with $d+3$ vertices. Using the expressions 
of $P(x)$ and $P^{+}(x)$ obtained in Theorem~\ref{theo:enum} and 
Theorem~\ref{theo:enum_oriented}, this yields the following corollary.  
\begin{corollary}
\label{cor:achiral}
Let $c^{-}(d+3,d)$ be the number of combinatorially different 
achiral $d$-polytopes with $d+3$ vertices. Then the generating 
function $P^{-}(x)=\sum_dc^{-}(d+3,d)x^{d+3}$ is equal to
{\small
$$
{\frac { \left( 2\,{x}^{11}+4\,{x}^{10}-2\,{x}^{9}-15\,{x}^{8}-5\,{x}^
{7}+23\,{x}^{6}+15\,{x}^{5}-17\,{x}^{4}-14\,{x}^{3}+4\,{x}^{2}+5\,x+1
 \right) {x}^{5}}{ \left( -1+x \right) ^{5} \left( 2\,{x}^{6}-4\,{x}^{
4}+4\,{x}^{2}-1 \right)  \left( x+1 \right) ^{3}}}.
$$
}
\end{corollary}
The first terms are $P^{-}(x)=x^5+7x^6+24x^8+62x^9+141x^{10}+287x^{11}+\ldots$.

\section{Complexity of the enumeration}
\label{sec:complexity}
The complexity model used here is the number of arithmetic operations, 
where an operation can be the addition of two integers 
of $\mathcal{O}(N)$ bits or the division of an integer of 
$\mathcal{O}(N)$ bits by a ``small'' integer, of $\mathcal{O}(\log(N))$ bits.

\begin{proposition}
\label{prop:complex}
The $N$ first coefficients counting combinatorial $d$-polytopes with $d+3$ 
vertices can be calculated in $\mathcal{O}(N\log(N))$ operations.

The $N$ first coefficients counting oriented $d$-polytopes with $d+3$ 
vertices can be calculated in $\mathcal{O}(N\log(N))$ operations.

The $N$ first coefficients counting combinatorial achiral $d$-polytopes with $d+3$ 
vertices can be calculated in $\mathcal{O}(N)$ operations.
\end{proposition}
\noindent\emph{Proof:}
In the proof we concentrate on the complexity of the extraction of the 
coefficients $c(d+3,d)$ from the expression of $P(x)$ given in 
Theorem~\ref{theo:enum} (the cases of $c^{+}(d+3,d)$ and $c^{-}(d+3,d)$ 
can be treated similarly). 

Given a generating function $f(x)$, 
we denote by $dev_N(f)$ the development of $f(x)$ up to power $x^N$. 
To calculate the $N$ first coefficients $c(d+3,d)$, it is sufficient to compute $dev_N\left((1-x)P(x)\right)$, where $P(x)=\sum_dc(d+3,d)x^{d+3}$. Indeed, the $d$th coefficient $f_d$ of $(1-x)P(x)$ verifies $f_{d+3}=c(d+3,d)-c(d+2,d-1)$. Hence, once the $N$ first coefficients $f_d$ are known, it takes $\mathcal{O}(N)$ operations to compute 
the $N$ first coefficients $c(d+3,d)$, that are calculated iteratively 
using 
$c(d+3,d)=c(d+2,d-1)+f_{d+3}$.

Multiplying the expression of $P(x)$ given in Theorem~\ref{theo:enum} by $(1-x)$ and 
then taking the $N$th truncation yields: 
\begin{eqnarray*}
dev_N\left((1-x)P(x)\right)&=&\sum_{e=1,e\ odd}^N \frac{\phi(e)}{4e}dev_N(G(x^e))+\sum_{e=1}^{N}\frac{\phi(e)}{2e}dev_N(H(x^e))\\
&&+dev_N(K(x)),
\end{eqnarray*}
where $G(x):=-\ln(1-2x^3/(1-2x)^2)$, 
$H(x):=-\ln(1-2x)+\ln(1-x)$ and $K(x)$ is an explicit rational 
function. 

As $K(x)$ is rational, its coefficients satisfy a linear recurrence 
with constant coefficients. Hence, finding $dev_N(K(x))$ requires 
$\mathcal{O}(N)$ operations. The development of $H(x)$ is explicit, 
$H(x)=\sum_n \frac{1}{n}(2^n-1)x^n$, so that finding $dev_N(H(x))$ 
requires $\mathcal{O}(N)$ operations. Let $G_n$ be the $n$th coefficient 
of $G(x)$. Then $xG'(x)=\sum_n nG_nx^n$. As $xG'(x)$ is a rational function, 
the coefficients $nG_n$ satisfy a linear recurrence with constant 
coefficients. Hence the calculation of $dev_N(G(x))$ requires 
$\mathcal{O}(N)$ operations. 

Once $dev_N(G(x))$, $dev_N(H(x))$ and $dev_N(K(x))$ are calculated, 
it remains to do the addition of $\phi(e)/(4e)dev_N(G(x^e))$ 
for odd $e$ from 1 to $N$. As $dev_N(G(x^e))$ has $\lfloor N/e\rfloor$ 
non zero coefficients, its addition takes $\mathcal{O}(N/e)$ operations. 
Hence the total cost of the addition is 
$\sum_{e=1}^N \mathcal{O}(N/e)=\mathcal{O}(N\ln(N))$ operations. 
Similarly the addition of $\phi(e)/(2e)dev_N(H(x^e))$ for $e$ 
from 1 to $N$ takes $\mathcal{O}(N\ln(N))$ operations. 
Finally it costs $\mathcal{O}(N\ln(N))$ operations to compute 
$dev_N\left((1-x)P(x)\right)$. 
$\Box$

\section{Asymptotic enumeration}
\label{sec:asym}
The asymptotic number of combinatorial $d$-polytopes, oriented $d$-polytopes and
achiral $d$-polytopes with $d+3$ vertices can be obtained from the
explicit formula of their generating functions, given respectively in
Theorem~\ref{theo:enum}, Theorem~\ref{theo:enum_oriented} and 
Corollary~\ref{cor:achiral}.

\begin{proposition}
The numbers $c(d+3,d)$ and $c^{+}(d+3,d)$ of combinatorial $d$-polytopes 
and oriented $d$-polytopes  with $d+3$ vertices have the 
asymptotic form:
$$
c(d+3,d)\sim\frac{\gamma^4}{4(\gamma-1)}\frac{\gamma^{d}}{d}\ \ \ \ \ \ \ c^{+}(d+3,d)\sim\frac{\gamma^4}{2(\gamma-1)}\frac{\gamma^{d}}{d}
$$  
where $\gamma^{-1}$ is the only real root of the equation $1-4x+4x^2-2x^3$, 
$\gamma\approx 2.839$. 

Let $\alpha$ be the unique positive root of $2x^6-4x^4+4x^2-1$ and $Q(x):=\frac{P^{-}(x)(1-4x^2+4x^4-2x^6)}{4x^5(3x^4-4x^2+2)}$. Then the number $c^{-}(d+3,d)$ of combinatorial achiral $d$-polytopes 
with $d+3$ vertices has the 
asymptotic form
$$
c^{-}(d+3,d)\sim \left(C+(-1)^dC'\right)\lambda^d
$$
where  $C:=Q(\alpha)\approx 12.1278$, $C':=Q(-\alpha)\approx 0.0346$ and $\lambda:=\alpha^{-1}\approx 1.6850$.
\end{proposition}
\noindent\emph{Proof:}
As for the proof of Proposition~\ref{prop:complex}, we only concentrate on 
the case of $c(d+3,d)$, (the proofs for $c^{+}(d+3,d)$ and $c^{-}(d+3,d)$ 
can be done with the same tools).  We use the framework of singularity 
analysis to derive an asymptotic estimate of the 
coefficients $c(d+3,d)$ from the development of $P(x)$ at its dominant 
singularity (for a generating function with nonnegative ceofficients, 
the dominant singularity is the least real value where $P(x)$ ceases to 
be analytic). From the expression of $P(x)$ given in 
Theorem~\ref{theo:enum}, it is easy to check that the dominant singularity 
of $P(x)$ is the real value $\rho$ such that $1-2\rho^3/(1-2\rho)^2=0$, 
i.e. the real root of $1-4x+4x^2-2x^3$. In addition, $P(z)$ has the 
following singular development at $\rho$, holding in a vicinity 
$\{|z-\rho|<\epsilon\}\cap\{z-\rho\notin \mathbf{R}^{+}\}$:
$$
P(z)\mathop{\sim}_{z\to\rho}\frac{1}{4(1-\rho)}\ln\left(\frac{1}{1-z/\rho}\right).
$$
Denote by $[x^n]f(x)$ the $n$th coefficient of a generating function 
$f(x)$. As the singular development of $P(z)$ holds in 
a ``Camembert''-neighbourhood of $\rho$ and as $P(z)$ is 
aperiodic\footnote{A generating function is aperiodic if it can be 
written as $x^ag(x^b)$ with $a\geq 0$ and $b\geq 2$ and $g(x)$ a 
generating function. It is clear that the presence of two consecutive 
positive coefficients is sufficient to be aperiodic. 
Hence $P(x)=x^5+7x^6+\ldots$ is aperiodic.}, transfer theorems of 
analytic combinatorics~\cite{flajolet} ensure that
$$
[x^n]P(x)\mathop{\sim}_{n\to \infty}[x^n]\frac{1}{4(1-\rho)}\ln\left(\frac{1}{1-x/\rho}\right)=\frac{1}{4(1-\rho)}\frac{\rho^{-n}}{n}.
$$
As $c(d+3,d)=[x^{d+3}]P(x)$, this yields the asymptotic estimate 
$c(d+3,d)\sim\frac{\gamma^{4}}{4(\gamma-1)}\frac{\gamma^d}{d}$, 
where $\gamma:=\rho^{-1}$. $\Box$
The following corollary follows directly from the fact that the growth 
constant of $c^{-}(d+3,d)$ is smaller than the growth constant of $c(d+3,d)$:
\begin{corollary}
 The quantity $c^{-}(d+3,d)/c(d+3,d)$, i.e.
 the probability  that a combinatorial $d$-polytopes with $d+3$ vertices 
is achiral, is asymptotically exponentially small.
\end{corollary}
For example, $c^{-}(d+3,d)/c(d+3,d)$ is less than $10\%$ for $d=10$
and less than $0.1\%$ for $d=20$.


\bibliographystyle{plain}
\bibliography{Polytopes}

\begin{thebibliography}{1}

\bibitem{flajolet}
P.~Flajolet and A.~Odlyzko.
\newblock Singularity analysis of generating functions.
\newblock {\em SIAM J. Discrete Math.}, 3:216--240, 1990.

\bibitem{grunbaum}
B.~Gr{\"u}nbaum.
\newblock {\em Convex {P}olytopes}, volume 221 of {\em Graduate Texts in Math.}
\newblock Springer-Verlag, New York, 2003.
\newblock Second edition prepared by V. Kaibel, V. Klee and G. M. Ziegler
  (original edition: Interscience, London 1967).

\bibitem{LLoyd}
K.~Lloyd.
\newblock The number of $d$-polytopes with $(d+3)$ vertices.
\newblock {\em Mathematika}, 17:120--132, 1970.

\end{thebibliography}

\end{document}